\documentclass[12pt]{amsart}
\usepackage{amsmath,amssymb,color,amscd}
\usepackage{xspace}

\usepackage{tikz}
\usetikzlibrary{positioning,arrows,cd}

\newcommand{\Var}{{\rm{Var}_{\mathbb{C}}}}

\newcommand{\chiorb}{{\chi^{\rm orb}}}

\def\1{\underline{1}}

\def\AA{{\mathbb A}}

\def\Z{{\mathbb Z}}

\def\RR{{\mathcal R}}
\def\KK{{\mathcal K}}
\def\AA{{\mathcal A}}

\def\GG{{\mathcal{G}}}
\def\RR{{\mathcal R}}

\newtheorem{theorem}{Theorem}

\newtheorem{lemma}{Lemma}
\newtheorem{proposition}{Proposition}
\newtheorem{definition}{Definition}

\newenvironment{corollary}
{\smallskip\noindent{\bf Corollary\/}.}{\smallskip\par}

\newenvironment{remark}
{\smallskip\noindent{\bf Remark\/}.}{\smallskip\par}

\title[On universal \& generalized orbifold Euler characteristics]{On the universal and generalized orbifold Euler characteristics}

\author{S.M.~Gusein-Zade}
\author{I.~Luengo}
\author{A.~Melle-Hern\'andez}
\author{A.~Viruel}
\thanks{The work of the S.G. (Sections~\ref{sec:orbifold_Euler} 
and~\ref{sec:reduction}) was supported by the grant 24-11-00124 of the Russian
Science Foundation. The work of I.L. and A.M. was
supported by Spanish grant
PID2020-114750GB-C32/AEI/10.13039/501100011033.
The work of A.V. was supported by AEI (Spain) grant PID2020-118753GB-I00/AEI/10.13039/501100011033, and by PAIDI 2020 (Andalusia) grant PROYEXCEL-00827.}


\begin{document}
\def\eps{\varepsilon}

\begin{abstract}
We discuss the universal orbifold Euler characteristic and generalized orbifold Euler characteristics corresponding to finitely generated groups $A$ (the {\em $A$-Euler characteristics}). We show that the collection of all $A$-Euler characteristics for $A$ of the form $A'\times\Z$ with finite $A'$ determine the universal orbifold Euler characteristic.
\end{abstract}

\maketitle

\section{Introduction}\label{sec:intro}
In all the paper, $\chi(\cdot)$ denotes {\em the additive Euler characteristic} defined as the alternating sum of the ranks of the cohomology groups with compact support:
\begin{equation}\label{eqn:Euler_definition}
\chi(X)=\sum_{q\ge 0}(-1)^q \cdot {\rm rk\,} H^q_c(X;\Z)\,.
\end{equation}

The (classical) orbifold Euler characteristic was defined by physicists in~\cite{Dixon_et_al1}
and~\cite{Dixon_et_al2} (and appeared in the mathematical literature in~\cite{AS} and~\cite{HH}).
There were defined some generalizations of this
notion. First, in~\cite{AS} and~\cite{Bryan_Fulman-1998}, there were defined,
so-called, higher order analogues of the orbifold
Euler characteristic. Then, in~\cite{Tamanoi} (see also~\cite{Farsi_Seaton-2011}, \cite{Farsi_Seaton-2022}), there were defined their generalizations corresponding to arbitrary
(or, rather, almost arbitrary) finitely generated groups. (The higher order orbifoid Euler characteristics correspond to free abelian groups.)
All of them were first defined for a manifold
(or a topological space) with an action of a finite group. The fact that they are really invariants
of of orbifolds is related with the, so-called,
induction property:~\cite{GLM-PEMS-2019}.
(A more simple proof valid in a more restrictive setting (sufficient for this paper) can be found
in~\cite[Theorem~1]{Simple_Tamanoi}.)

In~\cite{GLM-FAA-2018}, there was defined the,
so-called, universal orbifold Euler characteristic
$\chi^{\rm un}(\cdot)$ which takes values in a
certain ring $\RR$ freely generated, as an abelian 
group, by the isomorphism classes of finite groups.
It was shown that $\chi^{\rm un}$ is the universal
additive topological invariant of orbifolds.
This fact means, in particular, that the universal
orbifold Euler characteristic determines all the
generalized orbifold Euler characteristics corresponding to finitely generated groups.

There is a natural question, whether the generalized
orbifold Euler characteristics carry the same
information as the universal one. In other words,
whether the generalized orbifold Euler characteristics
corresponding to all finitely generated groups
(or to finitely generated groups considered in~\cite{Tamanoi}: see Section~\ref{sec:orbifold_Euler}) determine the
universal Euler characteristic.

In this paper we give the affirmative answer to
this question. For that, we reduce the problem to a 
purely algebraic one and give a solution for the latter.
The main statement of the paper is the following.

\begin{theorem}\label{theo:main}
 The collection of the generalized orbifold Euler characteristics $\chi^{(A\times\Z)}(X,G)$ corresponding to
 all the groups $A\times\Z$ with finite $A$ determines
 the universal orbifold Euler characteristic
 $\chi^{\rm un}(X,G)$.
\end{theorem}

This statement implies, in particular, that this collection
of generalized orbifold Euler characteristics determines
all other ones (corresponding to other finitely generated
groups).

\section{The universal orbifold Euler characteristic}\label{sec:universal}
The universal Euler characteristic of orbifolds (or, probably better to say, 
of spaces with actions of finite groups in the sense
of~\cite{GLM-PEMS-2019}) was defined in~\cite{GLM-FAA-2018}.

Let $G$ be a finite group. For a topological $G$-space $X$
(nice enough, say, homeomorphic to a locally closed
union of cells in a finite
$G$-CW-complex, see, e.\,g.,~\cite{TtD})
and for a subgroup $K\subset G$, let $X^K$ be the fixed
point set $\{x\in X: gx=x \text{ for all } g\in K\}$ of the subgroup $K$ and let $X^{(K)}$ be the set of points of $X$
whose isotropy subgroups $G_x=\{g\in G: gx=x\}$ coincide with $K$. Let ${Conjsub\,}G$ be the set of the conjugacy classes of subgroups of $G$. For a class $[K]\in {Conjsub\,}G$ (containing the subgroup $K$), let $X^{[K]}$ (respectively $X^{([K])}$) be the set of points of $X$
whose isotropy group contains a subgroup from the class $[K]$ (respectively, is a subgroup from $[K]$).

For a $G$-space $X$ and for a finite group $H\supset G$,
let ${\rm Ind}_G^H X$ be the induced $H$-space defined in the following way. Let us consider the Cartesian product
$H\times X$ and let us say that two points $(h_i,x_i)$,
$i=1,2$, of it are equivalent if there exists $g\in G$
such that $h_2=h_1g$, $x_2=g^{-1}x_1$. Then
${\rm Ind}_G^H X$ is the quotient of $H\times X$
by this equivalence relation with the natural $H$-action:
$h(h',x)=(hh', x)$. For finite groups $K\subset G\subset H$,
one has ${\rm Ind}_K^H={\rm Ind}_G^H\circ {\rm Ind}_K^G$ and ${\rm Ind}_G^H G/K=H/K$.

The universal Euler characteristic of orbifolds takes values
in the ring $\RR$ which can be defined as the Grothendieck ring of finite sets (``zero-dimensional varieties'') with actions of finite groups in the sense of~\cite{GLM-PEMS-2019}. This means the follwoing.
Two finite sets $X_1$ and $X_2$ with actions of finite groups
$G_1$ and $G_2$ respectively (that is a $G_1$-set
$(X_1,G_1)$ and a $G_2$-set
$(X_2,G_2)$) are called {\em isomorphic}
if there exists a bijection $h:X_1\to X_2$ and
a group isomorphism $\varphi:G_1\to G_2$ such that,
for $x\in X_1$ and $g\in G_1$, one has
$h(gx)=\varphi(g)h(x)$. Now $\RR$ is the abelian group
generated by the isomorphism classes $[(X,G)]$ of finite 
sets with actions of finite groups modulo the relations:
\begin{enumerate}
 \item[1)] if $Y$ is a $G$-invariant subset of a $G$-set $X$, then $[(X,G)]=[(Y,G)]+[(X\setminus Y,G)]$;
 \item[2)] if $G$ is a subgroup of a finite group $H$,
 then $[(X,G)]=[({\rm Int}_G^H X,H)]$.
\end{enumerate}
The multiplication in $\RR$ is defined by the Cartesian product:
$[(X_1,G_1)]\cdot[(X_2,G_2)]=[(X_1\times X_2,G_1\times G_2)]$.

The ring $\RR$ can be also described in the following way. Let $\GG$ be the set of all isomorphisms classes of finite groups. Then $\RR$ is 
the free Abelian group generated by the elements $T^{\KK}$ corresponding to the classes $\KK\in \GG$.
The element $T^{\KK}$ is represented by the one-point
$G$-set $G/G$ with $G\in\KK$: $T^{\KK}=[(G/G,G)]$. Thus an element of $\RR$
can be written as a finite sum of the form 
$\sum_{ {\KK}\in \GG} a_{\KK} T^{\KK}$, where $a_{\KK}\in \Z$.
The multiplication is defined by
$T ^{\KK'}\cdot T^{\KK''}=T^{\{G' \times G''\}}$
with $G'\in\KK'$, $G''\in\KK''$. 
According to the Krull--Schmidt theorem each finite group has a unique
 representation as the Cartesian product of indecomposable groups. 
 Let $\GG_{ind}$ be the set of  
 the isomorphisms  classes of indecomposable  finite groups.
The Krull--Schmidt theorem implies that $\RR$ is the polynomial ring $\Z [T^{\KK}]$ in the variables $T^{\KK}$ corresponding to 
the isomorphisms classes of the indecomposable finite groups.

\begin{definition}\label{def:universal-euler} (cf. \cite[Definition~7]{GLM-FAA-2018})
For a $G$-space $X$ ($G$ is a finite group),
the universal Euler characteristic is defined by:
\begin{eqnarray} \label{eq:universal-euler}
 \chi^{\rm un}(X,G)&=&\sum_{\KK\in \GG}
 \left(\sum_{{[K]\in {\rm Conjsub\,} G}
 \atop{[K]\in\KK}}
 \chi(X^{([K])}/G)\right)\cdot T^{\KK}\nonumber
 \\
 {\ }&=&
 \sum_{{[K]\in {\rm Conjsub\,} G}}
 \chi(X^{([K])}/G)\cdot T^{\{[K]\}},\label{eqn:universal-euler}
\end{eqnarray}
where $\{[K]\}\in\GG$ is the isomorphism class of the groups from $[K]$.
\end{definition}

The right-hand side of Equation~(\ref{eqn:universal-euler}) can
be written as an integral with respect to the Euler characteristic:
$$
 \chi^{\rm un}(X,G)=\int_{X/G}T^{\{G_x\}}d\chi\,,
$$
where $\{ G_x\}$ is the isomorphy class of the isotropy subgroup
of a point from the orbit $Gx$.

This definition can be extended to orbifolds (that is to spaces locally isomorphic to quotients $X/G$).
It is not difficult to see that $\chi^{\rm un}$ is an additive and multiplicative invariant of $V$-manifolds. As it was mentioned in
Introduction, $\chi^{\rm un}$ is a universal additive topological
invariant of orbifolds.

\section{Generalized orbifold Euler characteristics}\label{sec:orbifold_Euler}
The (classical) orbifold Euler characteristic was defined in~\cite{Dixon_et_al1} and~\cite{Dixon_et_al2} by the following two
equivalent equations:
$$
\chiorb(X,G)=\frac{1}{\vert G\vert}
 \sum_{{(g_0,g_1)\in G^2:}\atop{g_0g_1=g_1g_0}}\chi(X^{\langle g_0,g_1\rangle})
=\sum_{[g]\in {{\rm Conj\,}G}} \chi(X^{\langle g\rangle}/C_G(g))\,,
$$
where $g$ is a representative of the class $[g]$, $C_G(g)$ is the
centralizer $\{h\in G: hg=gh\}$ of the element $g$ in $G$, and $\langle \ldots\rangle$ is the subgroup of $G$ generated by the corresponding elements.
(A proof that these two equations are equivalent can be found in~\cite{HH}, see
also~\cite[Proposition~6]{RMS-2017}.)

For a $G$-space $X$, the {\em higher order orbifold Euler characteristics} were defined in~\cite{AS} and~\cite{Bryan_Fulman-1998} by the following equations
\begin{equation}\label{eqn:chi-k-orb}
 \chi^{(k)}(X,G)=
\frac{1}{\vert G\vert}\sum_{{{\bf g}\in G^{k+1}:}\atop{g_ig_j=g_jg_i}}\chi(X^{\langle {\bf g}\rangle})
=\sum_{[g]\in {{\rm Conj\,}G}} \chi^{(k-1)}(X^{\langle g\rangle}, C_G(g))\,,
\end{equation}
where $k$ is a non-negative integer (the order of the Euler characteristic), ${\bf g}=(g_0,g_1, \ldots, g_k)$
and (for the second, recurrent, definition)
$\chi^{(0)}(X,G)$ is given by the first definition
and coincides with the Euler characteristic
$\chi(X/G)$ of the quotient: see, e.\,g., \cite[Proposition~7]{RMS-2017}.
The usual orbifold Euler characteristic $\chi^{\rm orb}(\cdot,\cdot)$
is the orbifold Euler characteristic of order 1. (One may include the Euler-Satake characteristic $\chi^{ES}(\cdot)$ (see \cite{Satake}) into this series:
the characteristic $\chi^{ES}(X,G)$ can be interpreted in a certain way as $\chi^{(-1)}(X,G)$.)

For a finitely generated group $A$, 
the {\em (generalized) orbifold Euler characteristic}
corresponding to $A$ (or {\em the $A$-Euler characteristic}),
is defined by
\begin{equation}\label{eqn:chi(A)}
\chi^{(A)}(X,G)=
\frac{1}{\vert G\vert}\sum_{\varphi\in {\rm Hom\,}(A,G)}\chi(X^{{\varphi(A)}})\,,
\end{equation}
where the summation is over the set of homomorphisms
$\varphi\colon A\to G$, $X^{{\varphi(A)}}$ is the set
of point of $X$ fixed with respect to all the elements of the image
of $\varphi$ (\cite{Tamanoi}, see also~\cite{Farsi_Seaton-2011}). The value $\chi^{(A)}(X,G)$ is, in general, a rational number.
In~\cite{Tamanoi}, the author considered only the
generalized orbifold Euler characteristics
corresponding to finitely generated groups of the form $A=\Z\times A'$.
An advantage of these class of the
generalized orbifold Euler characteristics
is that they take values in $\Z$ (see Proposition~\ref{prop:A-Euler_for_product}).
Another advantage is the fact that there is a way to define
``motivic versions'' of them. This means that,
for a complex quasiprojective $G$-variety $X$,
one can define an element of the Grothendieck
ring $K_0(\Var)$ of complex quasiprojective varieties
whose Euler characteristic is equal to $\chi^{(A)}(X,G)$.

It is not difficult to see that, for the trivial group $G$ (i.\,e., $G=\{e\}$)
the $A$-Euler characteristic $\chi^{(A)}(X,\{e\})$
coincides with the usual Euler characteristic
$\chi(X)$. The definition directly implies
the additivity and the multiplicativity of the $A$-Euler characteristic
in the following sense: 
if $X_1$ is a $G$-invariant subset of a $G$ space $X$,
$X_2=X\setminus X_1$, then
$\chi^{(A)}(X,G)=\chi^{(A)}(X_1,G)+\chi^{(A)}(X_2,G)$; if $X_1$ is a $G_1$-space and $X_2$ is a $G_2$-space, then
\begin{equation}\label{eqn:multiplicat}
\chi^{(A)}(X_1\times X_2,G_1\times G_2)=\chi^{(A)}(X_1,G_1)\cdot\chi^{(A)}(X_2,G_2).
\end{equation}

\begin{proposition}\label{prop:A-Euler_for_product} (\cite[Proposition~2-1]{Tamanoi})
 If $A=A_1\times A_2$, then
 \begin{equation}\label{eqn;A-Euler_for_product}
  \chi^{(A)}(X,G)=\sum_{[\varphi]\in {\rm Hom}(A_1, G)/G}
  \chi^{(A_2)}\left(X^{\varphi(A_1)},C_G(\varphi(A_1))\right),
 \end{equation}
 where the group $G$ acts on the set ${\rm Hom}(A_1, G)$ by conjugation,
 $\varphi$ is a representative of the conjugacy class $[\varphi]$.
\end{proposition}

\begin{proof} (cf.\ the proofs of \cite[Proposition~2-1]{Tamanoi} and ~\cite[Proposition~6]{RMS-2017}) One has:
\begin{eqnarray*}
 \ &{\ }&\chi^{(A)}(X,G)=\frac{1}{\vert G\vert}
 \sum_{\varphi\in {\rm Hom}(A, G)}\chi(X^{ {\varphi(A)}})\\
 \ &=&\frac{1}{\vert G\vert}\sum_{[\varphi]\in{\rm Hom}(A_1,G)/G}
 \vert[\varphi]\vert\cdot\sum_{\psi\in{\rm Hom}\left(A_2,C_G(\varphi(A_1))\right)}
 \chi(X^{\varphi(A_1)\cdot\psi(A_2)})\\
 \ &=&\frac{1}{\vert G\vert}\sum_{[\varphi]\in{\rm Hom}(A_1,G)}
 \frac{\vert G\vert}{\vert C_G(\varphi(A_1))\vert}
 \sum_{\psi\in{\rm Hom}\left(A_2,C_G(\varphi(A_1))\right)}
 \chi((X^{\varphi(A_1)})^{\psi(A_2)})\\
 \ &=&\sum_{[\varphi]\in {\rm Hom}(A_1, G)/G}
  \chi^{(A_2)}(X^{\varphi(A_1)},C_G(\varphi(A_1)))\,.
\end{eqnarray*}
\end{proof}

\section{Reduction to a group theoretic problem}\label{sec:reduction}
One can consider the $A$-Euler characteristic as a function on the ring $\RR$, defining $\chi^{(A)}([X,G])$ ($X$ is a finite $G$-set) as
$\chi^{(A)}(X,G)$, where $X$ is considered as a zero-dimensional
topological space (with the discrete topology). The fact that this function is well-defined
is a consequence of the following {\em induction property}.
(It was proved first for the higher order orbifold Euler characteristics in~\cite[Theorem~1]{GLM-PEMS-2019}. See a more simple proof
for an arbitrary finitely generated group $A$ (in a purely topological
setting) in~\cite[Theorem~1]{Simple_Tamanoi}.)

\begin{proposition}
 Let $G\subset H$ be finite groups. For a $G$-space $X$ one has
 $$
 \chi^{(A)}(X,G)= \chi^{(A)}({\rm Ind}_G^H X,H).
 $$
\end{proposition}

From the definition of the universal Euler characteristic, it
follows that
$$
 \chi^{(A)}(X,G)= \chi^{(A)}(\chi^{\rm un}(X,G)).
 $$

Now let $\AA$ be a set of finitely generated groups (considered up to
isomorphisms). The fact that the generalized orbifold Euler characteristics $\chi^{(A)}$ for $A\in\AA$ define the universal
Euler characteristic $\chi^{\rm un}$ is equivalent to the fact
that their common kernel in $\RR$ is trivial, i.\,e., if 
$\chi^{(A)}(r)=0$ for all $A\in\AA$ ($r\in\RR$), then $r=0$. 

We shall show that this is the case for $\AA$ being the set of
groups $A$ of the form $A=A'\times \Z$ with finite $A'$. For that
we have to compute $\chi^{(A)}(r)$ for a finite sum
$r=\sum_{\KK\in\GG}a_{\KK}T^{\KK}$ ($a_{\KK}\in\Z$).
The following statement gives the answer for $A=A'\times\Z$.

\begin{proposition}
 One has
 $$
 \chi^{(A'\times\Z)}(G/G,G)=\left\vert{\rm Hom}(A',G)/G\right\vert\,.
 $$
\end{proposition}

\begin{proof}
 Proposition~\ref{prop:A-Euler_for_product} applied to $A_1=A'$, $A_2=\Z$ gives
 $$
 \chi^{(A'\times\Z)}(G/G,G)=\sum_{\varphi\in {\rm Hom(A',G)/G}}
 \chi^{(\Z)}((G/G)^{\varphi(A')}, C_G(\varphi(A')))\,
 $$
 On the other hand, $\chi^{(\Z)}(Y,H)=\chi(Y/H)$ and 
 $(G/G)^{\varphi(A')}$ is the one-point set. Therefore,
 $\chi^{(\Z)}((G/G)^{\varphi(A')}, C_G(\varphi(A')))=1$.
\end{proof}

This gives the following criterion for the $A$-Euler characteristics for $A$ of the form $A'\times\Z$ (with $A'$ finite)
to define the universal orbifold Euler characteristic:
if, for $a_{\KK}\in\Z$ such that finitely many of them are different from zero
and for any finite group $A'$, one has
$$
\sum_{\KK\in\GG}a_{\KK}\left\vert{\rm Hom}(A',G)/G\right\vert=0\,,
$$
($G$ is a representative of $\KK$), then all $a_{\KK}$ are equal to zero. This will be proved in Section~\ref{sec:proof}: Theorem~\ref{thm:linear_combination_vanish}.

\begin{remark}
 It is somewhat easier to prove that the $A$-orbifold Euler
 characteristics for all finitely generated groups $A$ determine
 the universal orbifold Euler characteristic. This can be deduced from
 the statement that, for any two non-isomorphic finite groups $G_1$ and $G_2$, there exists a finite group $A$ such that
 $\left\vert{\rm Hom}(A',G_1)/G_1\right\vert\ne
 \left\vert{\rm Hom}(A',G_2)/G_2\right\vert$ (a very particular case of Theorem~\ref{thm:linear_combination_vanish} below).
\end{remark}

\section{Counting morphisms between finite groups}\label{sec:proof}
We shall prove some statements in a somewhat more general setting than it is
necessary for the described aim of the paper.
In what follows, for every given group $G$, we chose a fixed subgroup $\omega(G)$ of $\operatorname{Aut}(G)$. The assignation $\omega$ is not required to fulfill any other extra hypothesis; it needs not to be functorial or consistent within any class of groups.

Now, if $A,B$ are groups, then there is a left-action of $\operatorname{Aut}(B)$ on the set $\operatorname{Hom}(A,B)$ given by post-composition, thus a left $\omega(B)$-action as well. 

\begin{definition}\label{def:representation}
Let $A,B$ be groups, then
\begin{enumerate}
    \item The set of $\omega$-representations of $A$ in $B$ is the set of $\omega(B)$-orbits:
    $$\operatorname{Rep}_\omega(A,B):=\omega(B)\backslash\operatorname{Hom}(A,B).$$

    \item The set of faithful $\omega$-representations of $A$ in $B$ is the set of $\omega(B)$-orbits:
    $$\operatorname{FRep}_\omega(A,B):=\omega(B)\backslash\operatorname{Mono}(A,B).$$
\end{enumerate}
\end{definition}

\begin{lemma}\label{lem:orbit_decomposition}
Let $A,B$ be groups, then
$$\operatorname{Rep}_\omega(A,B)=\bigsqcup_{K\trianglelefteq A} \operatorname{FRep}_\omega(A/K,B).$$
\end{lemma}
\begin{proof}
Notice that every group morphism $f\colon A\to B$ factors in a unique way as
$$
\begin{tikzcd}
A \arrow[dr,-{>>}, "p"] \arrow{rr}{f} && B, \\
{} & A/\ker(f) \arrow[ur,hook,"\overline{f}"]
\end{tikzcd}
$$
where $\overline{f}$ is a monomorphism. In other words, every group morphism $f\colon A\to B$ is completely determined by the pair $(\ker(f),\overline{f})$ where $\ker(f)\trianglelefteq A$ and $\overline{f}\colon A/\ker(f)\to B$ is a monomorphism. And given any pair $(K,\overline{f}),$ where $K\trianglelefteq A$ and $\overline{f}\colon A/K\to B$ is a monomorphism, there exists a unique group morphism $f\colon A\to B$ determined by the composition 
$$\begin{tikzcd}A\arrow[r,-{>>}, "p"]& A/K \arrow[r,hook, "\overline{f}"]& B.
\end{tikzcd}$$
Hence
\begin{equation}\label{eq:morphism_decomposition}
    \operatorname{Hom}(A,B)=\bigsqcup_{K\trianglelefteq A} \operatorname{Mono}(A/K,B).
\end{equation}

Finally, the $\omega(B)$-action on $\operatorname{Hom}(A,B)$ is via post-composition, and therefore, given a pair $(K,\overline{f})$ the $\omega(B)$-action only applies to $\overline{f}$. Indeed, if $\theta\in\omega(B)$ and $f\in \operatorname{Hom}(A,B)$ then $\ker(f)=\ker(\theta\circ f)$, hence the orbit $[f]\in \operatorname{Rep}_\omega(A,B)$ is fully determined by the pair $(\ker(f), [\overline{f}])$ where $\ker(f)$ does not depend on the representative of the $\omega(B)$-orbit $[f]$, and $[\overline{f}]\in \operatorname{FRep}_\omega(A/\ker(f),B).$ Therefore, the $\omega(B)$-action on $\operatorname{Hom}(A,B)$ is compatible with the decomposition given by Equation \eqref{eq:morphism_decomposition}, and the result follows. 
\end{proof}

A straightforward consequence of Lemma \ref{lem:orbit_decomposition} is the following:

\begin{corollary}\label{cor:size_orbit_decomposition}
Let $A,B$ be finite groups, then
$$|\operatorname{Rep}_\omega(A,B)|=\sum_{K\trianglelefteq A} |\operatorname{FRep}_\omega(A/K,B)|.$$    
\end{corollary}

\begin{remark}\label{rmk:size_orbit_decomposition}
Note that the assumption of finiteness for both $A$ and $B$ is solely to ensure that the sum on the right-hand side of the formula makes sense.  
\end{remark}

The following result shows that considering linear combinations involving the size of $\omega$-representation is equivalent to consider linear combinations involving the size of faithful $\omega$-representation.

\begin{lemma}\label{lem:reduce_linear_combinations}
Let $R$ be a ring, $a_1,\ldots,a_r\in R$, and $G_1,\ldots,G_r$ be finite groups. Then the following are equivalent:
\begin{enumerate}
    \item\label{lem:reduce_linear_combinations_1} For every finite group $A$, the following equation holds: $$\sum_{i=1}^s a_i|\operatorname{Rep}_\omega(A,G_i)|=0.$$

    \item\label{lem:reduce_linear_combinations_2} For every finite group $A$, the following equation holds: $$\sum_{i=1}^s a_i|\operatorname{FRep}_\omega(A,G_i)|=0.$$
\end{enumerate}
\end{lemma}
\begin{proof}
First assume \eqref{lem:reduce_linear_combinations_2}. Then, according to Lemma \ref{lem:orbit_decomposition}:
\begin{align*}
 \sum_{i=1}^s a_i|\operatorname{Rep}_\omega(A,G_i)| 
 &=\sum_{i=1}^s a_i\big(\sum_{K\trianglelefteq A} |\operatorname{FRep}_\omega(A/K,G_i)|\big)\\
 &=\sum_{K\trianglelefteq A}\big(\sum_{i=1}^s a_i|\operatorname{FRep}_\omega(A/K,G_i)|\big)\\
 &=0,
\end{align*}
hence \eqref{lem:reduce_linear_combinations_1} follows.

Now, we assume \eqref{lem:reduce_linear_combinations_1}, and prove that \eqref{lem:reduce_linear_combinations_2} follows by an inductive argument on $|A|$.

If $|A|=1$, then $A$ is the trivial group and $\operatorname{FRep}_\omega(A/K,G_i)=\operatorname{Rep}_\omega(A/K,G_i)$ for every $i$, hence there is nothing to prove.

Assume now that \eqref{lem:reduce_linear_combinations_1} holds for every finite group of size smaller that $n$, and let $A$ be a group with $|A|=n.$ Then, according to Lemma \ref{lem:orbit_decomposition}:
\begin{align*}
0&=&\sum_{i=1}^s a_i|\operatorname{Rep}_\omega(A,G_i)| 
= \sum_{i=1}^s a_i\big(\sum_{K\trianglelefteq A} |\operatorname{FRep}_\omega(A/K,G_i)|\big)\\
&=& \sum_{i=1}^s a_i|\operatorname{FRep}_\omega(A,G_i)| + \sum_{\{1\}\neq K\trianglelefteq A}\big(\sum_{i=1}^s a_i|\operatorname{FRep}_\omega(A/K,G_i)|\big),
\end{align*}
where if $\{1\}\neq K\trianglelefteq A$ then $|A/K|<n$ and $\sum_{i=1}^s a_i|\operatorname{FRep}_\omega(A/K,G_i)|=0$ by the induction hypothesis. Hence \eqref{lem:reduce_linear_combinations_2} holds for $A$ too.  
\end{proof}

Now our main result in this section.

\begin{theorem}\label{thm:linear_combination_vanish}
Let $R$ be a ring and $\{G_1,\ldots,G_s\}$ be a finite set of pairwise non-isomorphic finite groups such that $|\operatorname{FRep}_\omega(G_i,G_i)|$ is not a zero divisor in $R$, $i=1,\ldots,s.$ If $a_1,\ldots,a_s\in R$ are such that for every finite group $A$ the equation
$$\sum_{i=1}^s a_i|\operatorname{Rep}_\omega(A,G_i)|=0$$
holds, then $a_i=0$ for $i=1,\ldots,s.$
\end{theorem}
\begin{proof}
According to Lemma \ref{lem:reduce_linear_combinations}, we may assume that for every finite group $A$ the equation
\begin{equation}\label{eq:linear_combination_vanish}
\sum_{i=1}^s a_i|\operatorname{FRep}_\omega(A,G_i)|=0    
\end{equation}
holds.

We may also assume that the finite groups have been ordered so that 
$$|G_1|\leq |G_2|\leq\ldots\leq|G_s|,$$
and since the groups $G_i$ are pairwise non-isomorphic, this implies
$$
|\operatorname{FRep}_\omega(G_j,G_i)|=0 \text{\ \ for\ \  } j<i.
$$

Now, if $A=G_s$ then Equation \eqref{eq:linear_combination_vanish} gives rise to 
$$a_s|\operatorname{FRep}_\omega(G_s,G_s)|=0,$$
and since $|\operatorname{FRep}_\omega(G_s,G_s)|$ is not a zero divisor in $R$, this implies $a_s=0.$ So we can assume now that we are dealing with just $s-1$ pairwise non-isomorphic finite groups, and iterate this process to end that $0=a_s=a_{s-1}=\ldots=a_1.$
\end{proof}

Now Theorem~\ref{theo:main} follows from
Theorem~\ref{thm:linear_combination_vanish} applied to
$w(G)$ being the group of inner automorphisms of $G$ (and
$R=\Z$).

\begin{theorem}
 The collection of the generalized orbifold Euler characteristics $\chi^{(A\times\Z)}(X,G)$ corresponding to
 all the groups $A\times\Z$ with finite $A$ determines
 the universal orbifold Euler characteristic
 $\chi^{\rm un}(X,G)$.
\end{theorem}

This statement implies, in particular, that this collection
of generalized orbifold Euler characteristics determines
all other ones (corresponding to other finitely generated
groups).

\bigskip

\noindent S.M.~Gusein-Zade

\noindent Lomonosov Moscow State University, Faculty
of Mechanics and Mathematics and 
Moscow Center for Fundamental and Applied Mathematics,
GSP-1, Moscow, 119991, Russia

\noindent \&
National Research University ``Higher School of Economics'',
Usacheva street 6, Moscow, 119048, Russia. 

\noindent E-mail:
sabir\symbol{'100}mccme.ru

\medskip
\noindent I.~Luengo\\
Department of Algebra, Geometry and Topology,\\
Complutense University of Madrid, Madrid, Spain\\
E-mail: iluengo\symbol{'100}mat.ucm.es

\medskip
\noindent A.~Melle-Hern\'andez\\
Instituto de Matem\'atica Interdisciplinar (IMI),\\
Department of Algebra, Geometry and Topology,\\
Complutense University of Madrid, Madrid, Spain\\
E-mail: amelle\symbol{'100}mat.ucm.es

\medskip
\noindent A.~Viruel\\
Departamento de \'Algebra, Geometr\'ia y Topolog\'ia,\\
Universidad de M\'alaga, Campus de Teatinos,\\
29071--M\'alaga, Spain\\
Email: viruel@uma.es

\end{document}